\documentclass[12pt]{article}
\usepackage{hyperref}
\usepackage{amsfonts,latexsym,rawfonts,amsmath,amssymb,amsthm}
\usepackage{amsmath,amssymb,amsfonts,latexsym,lscape,rawfonts}
\usepackage[cm]{fullpage}
\usepackage{amsmath,amscd, float,times,rotating}
\usepackage{pb-diagram}
\usepackage{color}
\numberwithin{equation}{section}

\newcommand{\pd}[2]{\frac {\partial #1}{\partial #2}}

\newcommand{\al}{\alpha}
\newcommand{\bb}{\beta}
\newcommand{\la}{\lambda}
\newcommand{\La}{\Lambda}
\newcommand{\oo}{\omega}

\newcommand{\dd}{\delta}
\newcommand{\Na}{\nabla}

\def\ga{\gamma}

\newcommand{\ee}{\epsilon}
\newcommand{\si}{\sigma}

\newcommand{\beq}{\begin{equation}}
\newcommand{\eeq}{\end{equation}}
\newcommand{\beqs}{\begin{eqnarray*}}
\newcommand{\eeqs}{\end{eqnarray*}}
\newcommand{\beqn}{\begin{eqnarray}}
\newcommand{\eeqn}{\end{eqnarray}}
\newcommand{\beqa}{\begin{array}}
\newcommand{\eeqa}{\end{array}}

\def\td{\tilde}

\def\NN{{\mathbb N}}

\def\ri{\rightarrow}

\def\un{\underline}
\def\si{\sigma}
\def\pbp{\sqrt{-1}\partial\bar\partial}

\def\Del{{\Delta}}

\def\cA{{\mathcal A}}

\def\Vol{{\rm Vol}}
\newtheorem{prop}{Proposition}[section]
\newtheorem{theo}[prop]{Theorem}
\newtheorem{lem}[prop]{Lemma}

\newtheorem{rem}[prop]{Remark}

\newtheorem{defi}[prop]{Definition}

\newtheorem{q}[prop]{Question}

\title{{\bf\Large{Stability of K\"ahler-Ricci flow}}}

\author{Xiuxiong Chen$^{*}$ and  Haozhao Li$^\dagger$}

\date{}

\begin{document}
\bibliographystyle{plain}
\maketitle

 {\emph{Abstract}}  We prove the convergence of K\"ahler-Ricci
flow with some small initial curvature conditions. As applications,
we discuss the convergence of K\"ahler-Ricci flow when the complex
structure varies on a K\"ahler-Einstein manifold.\\



\section{Introduction}
\subsection{The motivation}

In \cite{[Ha82]}, R. Hamilton introduced the famous Ricci flow which
deforms any Riemannian metric in the direction of negative Ricci
curvature. If the Ricci flow exists globally, this will lead to the
existence of some canonical geometric structure (either Einstein
metrics or some non-trivial solitons). Unfortunately, the flow
usually will develop singularities in finite time. One key issue is
to study the formation of singularities  over finite time
\cite{[Ha95]}.  According to R. Hamilton, the singularity can be
divided roughly into  two types: the fast forming one (Type I) and
the slow forming one (Type II).  From the PDE point of view, the
Type I singularity is more ``gentle" where the curvature is
controlled by ${1 \over {T-t}}.\;$  We will focus our discussions on
the  ``fast forming" singularities
now.  \\

Let $g$ be any  Riemannian metric and we start the Ricci flow from
$g.\;$ Suppose the maximal existence time is $T(g) < \infty.\;$
Suppose further that this is a Type I singularity.   One intriguing
question one may ask is:  for any metric in a small neighborhood of
$g$, will the maximal existence time be $T(g)$ also? Will that be
Type I?   Note that Einstein metrics with positive scalar curvature,
or Ricci solitons (either stable or shrinking),  develop exactly
finite time Type I singularities.  The question above can be
simplified to

\begin{q}  \cite{[Chen]} Let $g$ be a stable or shrinking Ricci
soliton.  Does there exist a small neighborhood where Ricci flow
will still develop only Type I singularities? Will the Ricci flow
asymptotically converge to some Ricci soliton?
\end{q}

 This is a very interesting question.  Unfortunately, the
answer in full generality should be negative: Even if we assume a
metric $g$ is $C^\infty$ close to an Einstein metric with positive
scalar curvature, it is not clearly at all whether the flow will
only develop Type I singularities, and whether the flow will
converge to an Einstein metric after proper re-scaling. Nonetheless,
it is an intriguing question to find a set of suitable geometric
conditions so that this question has an affirmative answer. This
type of problem is called linear stabilities and there are lots of
research on this topic
(cf. \cite{[Na2]}, \cite{[DWW]}, \cite{[GIK]}, \cite{[Ye]} etc).  \\

One way to approach this problem is to associate certain functional
to Ricci flow.  In general setting, the first choice of course will
be Perelman's entropy functional. Unfortunately, Perelman's entropy
functional depends on the blowing up time $T$ and we have no idea
how $T$ will vary when
we vary initial metric.  \\

It is now clear why we should study this important problem on
K\"ahler manifolds first.  In any compact manifold, for the Ricci
flow initiating from any K\"ahler metric, the blowing up time is not
only a priori determined, but more importantly it can be explicitly
computed in terms of the first Chern class and the K\"ahler class.
Therefore, Question 1.1 will be highly interesting and feasible with
present technology when restricted to the K\"ahler setting.   A more
modest question one may ask is

\begin{q}  \cite{[Chen]}  Let $(M, g, J)$ be a K\"ahler-Einstein
manifold.  Does there exist a small neighborhood of the
K\"ahler-Einstein metric such that the K\"ahler-Ricci flow will
converge to a K\"ahler-Einstein metric?
\end{q}

If we don't perturb the complex structure, the answer is yes due to
the unpublished work of G. Perelman \cite{[Pere2]}, as well as
Tian-Zhu \cite{[TZ1]}. If we indeed perturb  the complex structure,
then the answer is not clear.  LeBrun-Simanca gives some results on
the related problem in \cite{[LS]}. An example of  G. Tian
\cite{[T97]} and S. K. Donaldson \cite{[Do2]} shows that the
Mukai-Umemura manifold $X$ admits a K\"ahler-Einstein metric, but
some small deformation of $X$ has no K\"ahler-Einstein metrics.
\\

The question we want to ask is of more general  than Question 1.2.

\begin{q}   \cite{[Chen]} Let $(M, g, J)$ be an ``almost K\"ahler-Einstein
manifold" in some natural sense. Will the K\"ahler-Ricci
flow converge to  a K\"ahler-Einstein metric (regardless of the
complex structure)?
\end{q}

This is the project we want to study in this paper. We will give
some stability results of K\"ahler-Ricci flow with respect to the
deformation of the underlying complex structures and  prove some new
convergence results with small energy conditions.\\

\subsection{On the deformation of complex structures}

Let $(M, [\omega])$ be a polarized compact K\"ahler manifold with
$[\omega]=2\pi c_1(M)>0 $ (the first Chern class) in this paper. One
of the main theorem we prove in this paper is:

\medskip

\begin{theo}\label{main1}Let $(M, g_{KE}, J_{KE})$ be a K\"ahler-Einstein  manifold with $c_1(M)>0,$ and no
non-zero holomorphic vector fields. For any K\"ahler metric
$\oo_g\in 2\pi c_1(M)\cap H^{1, 1}(M, J)$ with possibly different
complex structure $J$ satisfying \beq\|(g, J)-(g_{KE},
J_{KE})\|_{C^2}\leq \ee,\label{main1:eq}\eeq for sufficiently small
$\ee(g_{KE}, J_{KE})>0,$  the K\"ahler-Ricci flow with the initial
metric $(\oo_g, J)$ will converge exponentially fast to a
K\"ahler-Einstein metric.
\end{theo}

\begin{rem}Theorem \ref{main1} can be derived from a combination of N. Koiso's
results (cf. Proposition 10.1 in \cite{[K]}) and Perelman or
Tian-Zhu's results \cite{[TZ1]}. More precisely, N. Koiso proved
that if M has no non-zero holomorphic vector fields, for any
one-parameter complex deformation $J_t$ of complex structure
$J_{KE}$, there exists a sequence of Einstein metrics $g_t$ which
are  K\"ahler metrics compatible with $J_t$.\\
\end{rem}

Next we consider the case that $(M, J)$ has non-zero holomorphic
vector fields. As in \cite{[CLW]}, we need to assume $(M, J)$ is
pre-stable, which means that the complex structure doesn't jump
under the action of diffeomorphism group of $M$ (cf. Definition
\ref{prestable}). Under this assumption, we have the following
result:

\begin{theo}\label{main3}Let $(M, g_{KE}, J_{KE})$ be a K\"ahler-Einstein  manifold with $c_1(M)>0.$ For any K\"ahler metric
$\oo_g\in 2\pi c_1(M)\cap H^{1, 1}(M, J)$ with possibly different
complex structure $J$ satisfying the following conditions:
\begin{enumerate}
\item $(M, J)$ is pre-stable;
\item $([\oo_g], J)$ has vanishing Futaki invariant;
  \item $\|(g, J)-(g_{KE}, J_{KE})\|_{C^2}\leq \ee,$
for sufficiently small $\ee(g_{KE}, J_{KE})>0;$
\end{enumerate}
the K\"ahler-Ricci flow with the initial metric $(\oo_g, J)$ will
converge exponentially fast to a K\"ahler-Einstein metric.
\end{theo}

\begin{rem}Lebrun-Simanca proved the existence of
constant-scalar-curvature K\"ahler metrics for a deformation of
complex structures on a K\"ahler manifold (cf. Theorem 5 in
\cite{[LS]}). They assume that the Futaki invariant is
non-degenerate, which says that the linearization of the Futaki
invariant in the direction of the K\"ahler class is injective. Note
that the Futaki invariant of a K\"ahler-Einstein metric is never
non-degenerate(cf. \cite{[LS]}).  Our case seems to be complimentary
to the case considered in \cite{[LS]}.
\end{rem}

The proof of Theorem \ref{main1} and \ref{main3} follows directly
from Theorem \ref{main5} and \ref{main6} respectively(cf. Section 5,
6). The idea of the proof follows from our previous paper
\cite{[CLW]},
and we will discuss the details in Section \ref{sec1.3}. \\

\subsection{On the convergence of K\"ahler-Ricci flow}\label{sec1.3}

In our previous paper \cite{[CLW]}, we proved some convergence
theorems for the K\"ahler-Ricci flow with certain initial energy and
curvature conditions. Here we will refine those arguments and prove
the following type of stability results :

\begin{theo}\label{main5}Let $(M,  J)$ be a K\"ahler manifold with
$c_1(M)>0.$ For any $\ga, \La>0$, there exists a small positive
constant $\ee(\ga, \La)>0$ such that for any metric $g$ in the
subspace of K\"ahler metrics \beq\{\;\oo_g\in 2\pi
c_1(M)\;|\;\la_1(\oo_g)>1+\ga,\quad |Rm|(\oo_g)\leq \La, \quad
Ca(\oo_g)\leq \ee \}, \label{main5:eq}\eeq where $\la_1(g)$ is the
first eigenvalue of the metric $\oo_g$ and $Ca(\oo_g)$ denotes the
(normalized) Calabi energy, the K\"ahler-Ricci flow with the initial
metric $\oo_g$ will converge exponentially fast to a
K\"ahler-Einstein metric.\\
\end{theo}

\begin{rem}The assumption  (\ref{main5:eq})
implies that $(M, J)$ has no non-zero holomorphic vector fields.
Note that we don't assume the existence of K\"ahler-Einstein metrics
on $(M, J).\;$
\end{rem}

One might tempt to think that condition (\ref{main5:eq}) implies the
existence of K\"ahler-Einstein metrics. That might be true (one need
to address a potential collapsing issue with (\ref{main5:eq})),
except that it is a K\"ahler-Einstein metric with
 possibly different complex structures.   The difficulty is really
 created by the fact the space of complex structure modulo
 diffeomorphisms   is not a Hausdorff space. \\

 As in Theorem \ref{main3}, when $(M, J)$ has non-zero holomorphic vector fields, we need to
 assume the pre-stable condition.

\begin{theo}\label{main6}Let $(M,  J)$ be a K\"ahler manifold with
$c_1(M)>0.$  Suppose $(M, J)$ is pre-stable and the Futaki invariant
of the class $2\pi c_1(M)$ vanishes. For any $\La>0$, there exists
$\ee(\La)>0$ such that for any metric $g$ with its K\"ahler form
$\oo_g$ in the subspace of K\"ahler metrics \beq \{\;\oo_{g}\in 2\pi
c_1(M) \;  |\;|Rm|(\oo_g)\leq \La, \quad Ca(\oo_g)\leq \ee
\;\},\label{main6:eq}\eeq the K\"ahler Ricci flow with the initial
metric $\oo_g$ will converge exponentially fast to a
K\"ahler-Einstein metric.
\end{theo}

\begin{rem}Under the pre-stable condition on complex structures,
Phong-Sturm \cite{[PS1]} and Phong-Song-Sturm-Weinkove
\cite{[PSSW2]} proved some convergence results of K\"ahler-Ricci
flow with extra curvature conditions. We refer the readers to
\cite{[PSSW]} \cite{[N]}\cite{[FZ]} for more recent results on
K\"ahler-Ricci flow.
\end{rem}

This type of stability problems for the K\"ahler-Ricci flow was
initiated in \cite{[C1]} and later
 in \cite{[CLW]} with an assumption on the smallness of energy functional $E_1$
 or $E_0$. In this paper,  we replace this energy condition by the assumption that the Calabi energy is
sufficiently small. Unlike Theorem 1.5 in \cite{[CLW]}, we don't
need any conditions on the potential function of the initial
K\"ahler metric.\\

The proof of the main theorems is more tricky than that in
\cite{[CLW]}, but the ideas are the same. First by Sprouse's result
in \cite{[Sp]} the smallness of the Calabi energy implies that the
$L^{\infty}$ norm of the traceless Ricci curvature is small after a
short time(cf. Proposition \ref{sec7:prop}), which further implies
that the eigenvalue is strictly great than $1$ (cf. Lemma
\ref{lemeig}). The eigenvalue estimates can be used to prove the
exponential decay of the traceless Ricci curvature for a short time
(cf. Theorem \ref{theo:Ric1} and \ref{theo:Ric2}), which implies the
full curvature tensor is uniformly bounded for a short time by Yau's
estimates (cf. Theorem \ref{theobisec}). However, the boundedness of
the full curvature tensor in turn implies the rough curvature
estimates for the next time interval, and we can repeat the previous
arguments. Using this "iteration" idea we can actually prove that
the full curvature tensor is uniformly bounded for all time.\\

In our subsequent papers, we will remove the condition on the Futaki
invariant and  the bound of the full curvature tensor, and give more
general results on the relation between the pre-stable condition and
the convergence of K\"ahler-Ricci flow. Recently, Tian-Zhu
\cite{[TZ2]} proved very interesting and much stronger results on
stability of K\"ahler-Ricci flow by using Perelman's $W$-functional.\\

{\bf Acknowledgements}: The second named author would like to thank
Professor F. Pacard for the continuous support and encouragement
during the course of this work. The second named author would also
like to thank Professor W. Y. Ding and X. H. Zhu  for their help and
some enlightening discussions.
\medskip

\section{Preliminaries}
Let $M$ be a compact K\"ahler manifold with $c_1(M)>0.$ Choose an
initial K\"ahler metric $g$ with the K\"ahler form  $\oo_g\in 2\pi
c_1(M).\;$ By the Hodge theorem, any K\"ahler form in the same
K\"ahler class can be written as
$$\oo_{\varphi}=\oo_g+\pbp \varphi$$
for some real potential function $\varphi$ on $M$. The
K\"ahler-Ricci flow (cf. \cite{[Ha82]}) on a K\"ahler manifold $M$
is of the form
\begin{equation}
{{\partial g_{i \overline{j}}} \over {\partial t }} = - R_{i
\overline{j}}+g_{i\bar j},\qquad\forall\; i,\; j= 1,2,\cdots ,n.
\label{eq:kahlerricciflow}
\end{equation}
It follows that on the level of K\"ahler potentials, the
K\"ahler-Ricci flow becomes
\begin{equation}
{{\partial \varphi} \over {\partial t }} =  \log
{{\omega_{\varphi}}^n \over {\omega_g}^n } + \varphi - h_{g} ,\quad
\varphi(0)=0, \label{eq:flowpotential}
\end{equation}
where $h_{g}$ is defined by \beq {\rm Ric}({\omega_g})- \omega_g =
\sqrt{-1} \partial \overline{\partial} h_ {g} \quad {\rm
and}\quad\displaystyle \int_M\; h_{g} {\omega_g}^n = 0.
\label{norh}\eeq Let $\un R$ be the  the average of the scalar
curvature, which is a constant depending only on the K\"ahler class
and the underlying complex structure. Then the normalized Calabi
energy (cf. \cite{[Ca1]}\cite{[Ca2]}) is defined by \beq
Ca(\oo_{\varphi})=\frac 1V\int_M\; (R(\oo_{\varphi})-\un
R)^2\;\oo_{\varphi}^n.\eeq Since the K\"ahler metric $\oo_{\varphi}$
is in the canonical class, we can check that
$$Ca(\oo_{\varphi})=\frac 1V\int_M\; |Ric(\oo_{\varphi})-\oo_{\varphi}|^2\;\oo_{\varphi}^n.$$
Define the Futaki invariant by
$$f_M([\oo_g], X)=\int_M\; X(h_g)\;\oo_g^n,$$
for any holomorphic vector field $X$ on $M$. It is well-known that
the Futaki invariant doesn't depend on the particular representative
we choose in the K\"ahler class.\\

 We recall some basic results from our previous papers.
 First, by the tensor maximum principle we have the following basic
lemma:

\begin{lem}\label{lemRm}(\cite{[C1]}\cite{[CLW]})
Suppose that the curvature of the initial metric satisfies the
following condition
$$\left \{\begin{array}{lll}|Rm|(0)&\leq& \La,\\
|Ric-\oo|(0)&\leq&\ee.
\end{array}\right. $$
there exists a constant $T(\La)>0$, such that we have the following
bound for the evolving K\"ahler metric $g(t)(0\leq t\leq 6T)$
\beq \left \{\begin{array}{lll}|Rm|(t)&\leq& 2\La,\\
|Ric-\oo|(t)&\leq& 2\ee.
\end{array}\right. \label{1}\eeq
\end{lem}

Lemma \ref{lemRm} is slightly different from Lemma 12 in
\cite{[C1]}, but the idea of the proof is the same.  We remind the
readers that $T$ doesn't depend on the bound of the traceless Ricci
curvature, which is useful for the proofs of the
main theorems.\\

Now we  state a parabolic version of Moser iteration argument (cf.
\cite{[CT1]}).

\begin{theo}\label{lemmoser} Suppose the Poincare constant and
the Sobolev constant of the evolving K\"ahler metrics $g(t)$ are both
uniformly bounded by $\sigma$, and the scalar curvature $R(g(t))$
has a uniform lower bound $\La$. If a nonnegative function $u$
satisfying the following inequality
$$\pd {}tu\leq \Del u+f(t, x)u, \;\; \forall a<t<b,$$
where $|f|_{L^p(M, g(t))}$ is uniformly bounded by some constant $c$
for some $p>\frac n2$, then for any $t\in (a, b)$ and $\tau\in (0,
b-a)$, we have\footnote{The constant $C$ can be different at
different places with possibly some lower indices. The notation
$C(A, B, ...)$ means that the constant $C$ depends only on $A, B,
...$.}
$$u(t)\leq \frac {C(n, \sigma, c, \La)}{\tau^{\frac {m+2}{4}}}
\Big(\int_{t-\tau}^t\int_M \;u^2 \oo_{\varphi}^n\wedge
ds\Big)^{\frac 12}.$$
\end{theo}

\begin{rem}Recently, Q. Zhang \cite{[Z1]}and R. Ye \cite{[Y1]} proved that the Sobolev
constant is uniformly bounded along the K\"ahler-Ricci flow without
any assumptions. Here we don't need to use this result.
\end{rem}

The parabolic Moser iteration theorem is the main tool in the paper.
It can be applied to control the pointwise norm of the traceless
Ricci curvature, provided that the Calabi energy is sufficiently
small. Once we have the bound of the traceless Ricci curvature, we
can estimate the potential function  and the full curvature tensor
along the K\"ahler-Ricci flow by Yau's estimates:

\begin{theo}\label{theoRm}(cf. \cite{[CLW]}, \cite{[CT2]}, \cite{[Yau]}) For any positive constants $\La, B>0$ and small $\eta>0$, there exists a constant $C_{1}$ depending
only on $\La, B, \eta$ such that if the {background metric} $\oo$
satisfies
$$|Rm|(\oo)\leq \La, \qquad |Ric(\oo)-\oo|\leq \eta, $$
and the potential function $|\varphi(t)|, |\dot\varphi(t)|\leq B,$
then
$$|Rm|(t)\leq C_{1}(B, \La, \eta).$$
\end{theo}

\medskip

 We state the following well-known result on the estimate of
the Sobolev constant. The readers are referred to
\cite{[Cro]}\cite{[LY]} for details.

\begin{theo}\label{theo:Sob}Let $(M, g)$ be a compact m-dimensional Riemannian manifold.
Suppose  $Ric\geq -\La, \Vol(M)\geq \nu>0,$  and the diameter
$diam(M)\leq D$, then there exists a constant $C_S(\La, \nu, D)>0$
such that for any function $f\in C^{\infty}(M),$ we have
$$\Big(\int_M\; |f|^{\frac {2m}{m-2}}\;dV_g\Big)^{\frac {m-2}{m}}\leq C_S\Big(\int_M\; |\Na f|^2\;dV_g+\int_M\;|f|^2\;dV_g\Big).$$

\end{theo}

\section{Estimates}
In this section, we will prove several results which will be useful
in the proof of the main theorems.

\subsection{The first eigenvalue of the Laplacian
operator}\label{sec:eigen}

To prove the exponential decay of the traceless Ricci curvature, we
need to estimate the first eigenvalue of the evolving Laplacian
operator. The calculation of the eigenvalue along the Ricci flow is
well-known in literatures(cf. \cite{[CaoX]} for example).

\begin{lem}\label{lemeig}Let $\la_1(t)$ be the first eigenvalue of the Laplacian operator acting on
functions with respect to the metric $\oo_{\varphi}$ along the
K\"ahler-Ricci flow.
\begin{enumerate}
  \item If $|Ric-\oo|(t)\leq \ee$ for $t\in [0, T]$, then $$\la_1(t)\geq
\la_1(0)e^{-3n\ee t},\quad \forall t\in [0, T].$$
  \item If $|Ric-\oo|(t)\leq \ee e^{-\al t}$ for some $\al>0$ and for all $t\in [0, T]$, then $$\la_1(t)\geq
\la_1(0)e^{-\frac {3n}{\al}\ee(1-e^{-\al t})},\quad \forall t\in [0,
T].$$
\end{enumerate}

\end{lem}
\begin{proof}
Let $\la_1(t)$ be a eigenvalue of $\Delta_{\varphi}$ with
$-\Delta_{\varphi} f(t)=\la_1(t)f(t),$ where $f(t)$ is a smooth
function  satisfying the normalization condition
$$\int_M\; f(t)^2\oo_{\varphi}^n=1.$$
Taking the derivative with respect to $t$, we have \beq
\int_M\;\Big( 2f\pd ft+f^2\Delta_{\varphi} {\pd {\varphi}t}\Big)
\,\oo_{\varphi}^n=0.\label{eq: eigen1}\eeq Observe that
$$\la_1(t)=\int_M\; |\Na f|^2\oo_{\varphi}^n,$$ we calculate the
derivative of $\la_1(t)$ \beqs \frac {d\la_1(t)}{dt}&=&-\frac
{d}{dt}\int_M\;f\Delta_{\varphi}f \;\oo_{\varphi}^n\\
&=&\int_M\; \Big(-\pd ft\Delta_{\varphi }f-f\pd {}t(\Delta_{\varphi}
f)-f\Delta_{\varphi}f \Delta_{\varphi} \pd
{\varphi}t\Big)\;\oo_{\varphi}^n\\
&=& \int_M\; \Big(-\pd ft\Delta_{\varphi}f-f\Delta_{\varphi} \pd
{f}t +f\Big(\pd {\varphi}t\Big)_{i\bar j}f_{j\bar
i}-f\Delta_{\varphi}f
\Delta_{\varphi} \pd {\varphi}t\Big)\;\oo_{\varphi}^n\\
&=&\la_1\int_M\;\Big(2 f\pd ft+f^2 \Delta_{\varphi} \pd
{\varphi}t\Big)\;\oo_{\varphi}^n+\int_M\; f\Big(\pd
{\varphi}t\Big)_{i\bar j}f_{j\bar i}\;\oo_{\varphi}^n. \eeqs
Applying (\ref{eq: eigen1}), we have that \beqn \frac
{d\la_1(t)}{dt}&=&\int_M\; f\Big(\pd {\varphi}t\Big)_{i\bar
j}f_{j\bar i}\oo_{\varphi}^n\nonumber\\&=&\int_M\;\Big(
(Ric(\oo_{\varphi})-\oo_{\varphi})(\Na f, \Na f)+fR_{, \bar \bb}f_{\bb}\Big)\;\oo_{\varphi}^n\nonumber\\
&=& \int_M\; (Ric(\oo_{\varphi})-\oo_{\varphi})(\Na f, \Na
f)\;\oo_{\varphi}^n-\int_M\; R|\Na f|^2
\oo_{\varphi}^n+\la_1(t)\int_M\;R f^2\oo_{\varphi}^n
.\label{eq:lemeig1}\eeqn The assumption (1) implies
$$Ric(\oo_{\varphi})-\oo_{\varphi}\geq -\ee\;\oo_{\varphi},\quad n-n\ee\leq R\leq n+n\ee.$$
Hence, we get the inequalities
 \beqs \frac {d\la_1(t)}{dt}
&\geq& -\ee \la_1(t)-(n+n\ee)\la_1(t)+(n-n\ee)\la_1(t)\\
&\geq&-3n\ee\la_1(t).\eeqs The first part of the lemma follows
immediately. Similarly we can prove the second part.
\end{proof}

\subsection{The pre-stable condition}
In this section, we estimate the first eigenvalue when $M$ has
non-zero holomorphic vector fields. Here we follow closely the
argument in \cite{[CLW]}. First, we recall the following definition
in \cite{[CLW]}:
\begin{defi}\label{prestable} The complex structure $J$ of $M$
is called pre-stable, if no complex structure in the closure of its
orbit of diffeomorphism group contains larger (reduced) holomorphic
automorphism group.

\end{defi}
\medskip
\begin{rem}
The "pre-stable" condition was defined in \cite{[CLW]}, but it is
well-known in previous literatures. In the statement of Theorem 1.8
of \cite{[C1]}, the first named author used this condition to study
the convergence to K\"ahler-Ricci flow, and in \cite{[PS1]}
Phong-Sturm defined a stability condition as "condition (B)". These
two definitions are essentially the same and we called it
"pre-stable" in \cite{[CLW]}.
\end{rem}

Now we recall some basic facts of the first and second eigenvalues
of the Laplacian operator acting on functions. For any smooth
function $f\in C^{\infty}(M)$ we have \beq \int_M\; |\Na \Delta_g
f|^2 dV_g\geq \la_1\;\int_M\; |\Delta_g f|^2 dV_g, \quad \forall
f\in C^{\infty}(M) \label{sec3.2:eq:1},\eeq where $\la_1$ is the
first eigenvalue of $\Delta_g$. If we assume the function $f$ is
perpendicular to the first eigenspace of $\Delta_g$ then we get \beq
\int_M\; |\Na \Delta_g f|^2 dV_g\geq \la_2\;\int_M\; |\Delta_g f|^2
dV_g. \label{sec3.2:eq:2}\eeq These inequalities can be proved by
the eigenvalue decomposition of the function $f$. Note that for a
K\"ahler-Einstein manifold $(M, \oo_{KE})$ with non-zero holomorphic
vector fields, it is well-known that the first eigenvalue $\la_1=1$
and the first eigenspace, which we denote by $\eta(M)$, is
isomorphic to the space
of holomorphic vector fields.\\

Now we need the following convergence result of a sequence of
K\"ahler metrics, which is well-known in literature (cf.
\cite{[PS1]}, \cite{[T1]}).

\begin{prop}\label{compactness} Let $M$ be a compact K\"ahler manifold. Let
$(g(t), J(t))$ be any sequence of metrics $g(t)$ and complex
structures $J(t)$ such that $g(t)$ is K\"ahler with respect to
$J(t)$. Suppose  the following is true:
\begin{enumerate}\item For some integer $k\geq 1$, $|\Na^lRm|_{g(t)}$ is uniformly bounded for
any integer $l (0\leq l< k)$;

\item The injectivity radii $i(M,
g(t))$ are all bounded from below;

\item There exist two
uniform constant $c_1$ and $c_2$ such that $0<c_1\leq \Vol(M,
g(t))\leq c_2$.
\end{enumerate}
Then there exists a subsequence of $t_j$, and a sequence of
diffeomorphism $F_j: M\ri M$ such that the pull-back metrics $\tilde
g(t_j)=F_j^*g(t_j)$  converge in $C^{k, \al}(\forall \,\al\in (0,
1))$ to a  $C^{k, \al}$ metric $g_{\infty}$. The pull-back complex
structure tensors $\tilde J(t_j)=F_j^*J(t_j)$ converge in $C^{k,
\al}$ to an integral complex structure tensor $\tilde J_{\infty}$.
Furthermore, the metric $ g_{\infty}$ is K\"ahler with respect to
the complex structure $\tilde J_{\infty}$.

\end{prop}

These being understood, we have the:

\begin{theo}\label{theo:prestable}Suppose that $(M, J)$ is
pre-stable. For any $\La_0, \La_1>0$,   there exists $\ee>0$
depending only on $\La_0$ and $\La_1$ such that for any metric
$\oo\in 2\pi c_1(M),$ if \beq|Ric(\oo)-\oo|\leq \ee,\;\;
|Rm|(\oo)\leq \La_0,\;\;|\Na Rm|(\oo)\leq \La_1, \label{r1}\eeq then
for any smooth function $f$ satisfying \beq \int_M\; f\oo^n=0 {\;\;{
and}\;\;} \int_M\; X(f)\oo^n=0, \qquad \forall X\in \eta(M,
J),\label{theo:pre eq3}\eeq we have the following \beqn \int_M\;
|\Na f|^2\oo^n&>&(1+\ga(\ee, \La_0, \La_1))\int_M\;
|f|^2\oo^n,\label{theo:pre eq1}\\\int_M\; |\Na\Delta f|^2\oo^n
&>&(1+\ga(\ee, \La_0, \La_1))\int_M\; |\Delta
 f|^2\oo^n,\label{theo:pre eq2}\eeqn
where $\ga>0$ depends only on $\ee, \La_0$ and $\La_1.$

\end{theo}
\begin{proof}The inequality (\ref{theo:pre eq1}) was proved in \cite{[CLW]}, and  the argument also works for (\ref{theo:pre eq2}).
For the readers' convenience, we give the details here.

 Suppose not, for any positive numbers $\ee_m\ri 0$, there exists a sequence
of K\"ahler metrics $\oo_m\in 2\pi c_1(M)$ such that \beq
|Ric(\oo_m)-\oo_{m}|\leq \ee_m,\;\; |Rm|(\oo_m)\leq
\La_0,\;\;\;|\Na_m Rm|(\oo_m)\leq \La_1,\label{r2}\eeq where the
smooth functions $f_m$ satisfy \beq \int_M\; f_m\oo_m^n=0,  \quad
\int_M\; X(f_m)\oo_m^n=0, \qquad \forall X\in \eta(M,
J),\label{a3}\eeq \beq\int_M\; |\Na_m \Delta_m
f_m|^2\oo_m^n<(1+\ga_m)\int_M\; |\Delta_m
f_m|^2\oo_m^n,\label{eq5.25}\eeq where $0<\ga_m\ri 0.$ Without loss
of generality, we may assume that \beq \int_M\; |\Delta_m f_m|^2
\omega_m^n = 1, \qquad \forall m \in \NN, \label{a1}\eeq which means
\beq \int_M\; |\Na_m \Delta_m f_m|^2\oo_m^n\leq  1 + \gamma_m <
2.\label{a2} \eeq Then, $f_m $ will converge weakly in $W^{3,2}$ if
$(M, \omega_m)$ converges. However, according to Proposition
\ref{compactness}, $(M, \omega_m, J)$ will converge in $C^{2,
\al}(\al\in (0, 1))$ to $(M, \omega_\infty, J_\infty).\;$  In fact,
by (\ref{r2}) the diameters of $\oo_m$ are uniformly bounded. Note
that all the metrics $\oo_m$ are in the same K\"ahler class, the
volume is fixed. Then by (\ref{r2}) again, the injectivity radii are
uniformly bounded from below. Therefore, all the  conditions of
Proposition \ref{compactness} are satisfied.

 Note that the complex structure
$J_\infty$ lies in the closure of the orbit of diffeomorphisms,
while $\omega_\infty$ is a K\"ahler-Einstein metric of $(M,
J_\infty)$.   By the standard deformation theorem in complex
structures, we have
\[
   \dim  Aut_r (M, J) \leq \dim Aut_r(M, J_\infty).
\]
By abusing notation, we can write
\[
Aut_r(M, J) \subset Aut_r(M, J_\infty).
\]
By our assumption of pre-stable of $(M, J)$, we have the inequality
the other way around. Thus, we have
\[
\dim Aut_r(M, J) = \dim Aut_r(M, J_\infty),\;\;\;{\rm or}\;\;\;
Aut_r(M, J) = Aut_r(M, J_\infty).
\]
Now, let $f_\infty $ be the $W^{3,2}$ limit of $f_m$, then by
(\ref{a1}) and (\ref{a2}) we have
\[
  1 \leq  |f_\infty|_{W^{3,2}(M,\, \omega_\infty)} \leq C,
  \]
  and by (\ref{a3}) we have
  \[
  \int_M f_\infty \omega_\infty^n = 0, \qquad \int_M\; X(f_\infty) \omega_\infty^n =
  0,\qquad \forall X\in \eta(M, J_{\infty}).
  \]
  Thus, $f_\infty$ is a non-trivial function.  Since $\omega_\infty$ is a K\"ahler-Einstein metric, we have
  \[
  \int_M\; \theta_X f_\infty \omega_\infty^n = 0,
  \]
  where \[
  {\cal L}_X \omega_\infty =\pbp\theta_X.
  \]
  This implies that $f_\infty $ is perpendicular to the first eigenspace of $\triangle_{\omega_\infty}.\;$
  In other words, there is a $\delta > 0$ such that
  \[
     \int_M |\nabla_{\infty} \Delta_{\infty}f_\infty|^2 \omega_\infty^n > (1+ \delta) \int_M |\Delta_{\infty}f_\infty|^2 \omega_\infty^n = 1+ \delta.
  \]
 However, this contradicts the following fact:
 \beqs  \int_M \;|\nabla \Delta_{\infty}f_\infty|^2 \omega_\infty^n  & \leq & \displaystyle \lim_{m\rightarrow \infty}
  \int_M |\Na\Delta_m f_m|^2 \omega_m^n \\
 & \leq &  \displaystyle \lim_{m\rightarrow \infty} (1+ \ga_m) \int_M |\Delta_m f_m|^2 \omega_m^n = 1.
 \eeqs
The lemma is then proved.
\end{proof}

\def\vp{\pd {\varphi}t}

\subsection{The exponential decay of traceless Ricci curvature}
In this section, we will use the estimates of the first eigenvalue
in the previous subsection to prove the exponential decay of
traceless Ricci curvature.

\begin{theo}\label{theo:Ric1}Suppose for any time  $t\in [0, T]$,
$$|Rm|(t)\leq \La, \quad |Ric-\oo|(t)\leq H\ee,\quad \la_1(t)\geq 1+\gamma>1.$$
 Then for fixed $\tau<T$, there exists
$\ee_0(H, \ga)>0$ such that if $\ee\in (0, \ee_0)$ we have
$$|Ric-\oo|(t)\leq C_2(\La, \ga, \tau)H\ee e^{-\frac {\ga}4 t}, \quad t\in [\tau, T].$$

\end{theo}

\medskip

\begin{proof}
First, we prove that the Calabi energy decays exponentially \beq
Ca(t)\leq e^{-\ga t }Ca(0),\quad \forall t\in [0,
T].\label{sec3.3:eq1}\eeq In fact, direct calculation shows \beqn
&&\frac d{dt}Ca(t)\nonumber\\&=&\frac 1V\int_M\;\Big(
2\Delta_{\varphi} \vp\pd
{}t\Big(\Delta_{\varphi} \vp\Big)+\Big(\Delta_{\varphi} \vp\Big)^3\Big)\;\oo_{\varphi}^n\nonumber\\
&=&\frac 1V\int_M\; \Big(-2\Delta_{\varphi}\vp \Big|\Na\bar \Na
\vp\Big|^2-2\Big|\Na\Delta_{\varphi}
\vp\Big|^2+2\Big|\Delta_{\varphi} \vp\Big|^2+\Big(\Delta_{\varphi}
\vp\Big)^3\Big)\;\oo_{\varphi}^n.\label{eq:Ric1:1}
 \eeqn
Combining (\ref{sec3.2:eq:1}) with the assumptions, we have \beq
\frac d{dt}Ca(t)\leq -(2\ga-3nH \ee)Ca(t).\label{eq:Calabi}\eeq
If we choose $\ee$ sufficiently small, the inequality (\ref{sec3.3:eq1}) follows immediately.\\

Now applying the parabolic Moser iteration, the estimate of the
Calabi energy implies the exponential decay of the traceless Ricci
curvature. In fact, since the traceless Ricci tensor
$u=|Ric-\oo|^2(t)$ satisfies the following inequality:
$$\pd ut\leq \Delta u+c(n)|Rm|u,$$
By the parabolic Moser iteration Theorem \ref{lemmoser}, for $t\in
[\tau, T]$ we have \beqs |Ric-\oo|^2(t)&\leq&\frac
{C(\La)}{\tau^{\frac {m+4}{2}}} (\int_{t-\tau}^t\int_M \;|Ric-\oo|^4
\oo_{\varphi}^n \wedge ds)^{\frac
12}\\
&\leq &\frac {H\ee C(\La)}{\tau^{\frac {m+4}{2}}}
\Big(\int_{t-\tau}^t\int_M \;|Ric-\oo|^2 \oo_{\varphi}^n\wedge
ds\Big)^{\frac
12}\\
&=&\frac {H\ee C(\La)}{\tau^{\frac {m+4}{2}}}
\Big(\int_{t-\tau}^t\int_M \;\Big(\Delta_{\varphi} \vp\Big)^2
\oo_{\varphi}^n\wedge ds\Big)^{\frac 12}
\\&\leq&\frac {C(\La)H\sqrt{Ca(0)}\;\ee}{\sqrt{\ga}\tau^{\frac {m+4}{2}}}e^{-\frac
{\ga}2(t-\tau)}. \eeqs  Here we have used the fact that the Sobolev
constant is uniformly bounded since the traceless Ricci curvature is
small by the assumption. Note that at the initial time $Ca(0)\leq
H^2\ee^2,$ the above inequality implies
$$|Ric-\oo|(t)\leq \frac {C(\La)H\ee}{\ga^{\frac 14}\tau^{\frac {m+4}{4}}}e^{-\frac
{\ga}4(t-\tau)}.$$ The theorem is proved.
\end{proof}

When $M$ has non-zero holomorphic vector fields, we need to use the
pre-stable condition to get the exponential decay of the Ricci
curvature.
\begin{theo}\label{theo:Ric2}Suppose that $M$ is pre-stable and the Futaki invariant vanishes, and for any time  $t\in [0, T]$,
$$|Rm|(t)\leq \La, \quad |Ric-\oo|(t)\leq H\ee.$$
Then for fixed $\tau<T$, there exists $\ee_0(H, \ga)>0$ such that if
$\ee\in (0, \ee_0)$ we have
$$|Ric-\oo|(t)\leq C_3(\La, \ga, \tau)H\ee e^{-\frac {\ga}4 t}, \quad t\in [\tau, T].$$
Here $\ga$ is the constant obtained in Theorem \ref{theo:prestable}.

\end{theo}
\begin{proof}The argument is essentially the same as the proof of Theorem
\ref{theo:Ric1}. In fact, by Shi's results in \cite{[Shi]} all the
derivatives of the curvature tensor are bounded after a short time
$t=\frac {\tau}2,$ where $\tau$ is a fixed number as in Theorem
\ref{theo:Ric1}. Therefore the assumption (\ref{r1}) in Theorem
\ref{theo:prestable} is satisfied for $t\in [\frac {\tau}2, T]$ and
we have the inequalities (\ref{theo:pre eq1}) and (\ref{theo:pre
eq2}). Note that the function $\pd {\varphi}t-c(t)(c(t)=\frac
1V\int_M\; \pd {\varphi}t\;\oo_{\varphi}^n)$ satisfies the condition
(\ref{theo:pre eq3}), since the Futaki invariant vanishes
$$ f_M(X)=\int_M\; X\Big(\pd {\varphi}t\Big)=0,\quad \forall\;X\in \eta(M).$$
Therefore, we have the inequality
$$\int_M\; \Big|\Na\Delta_{\varphi} \pd {\varphi}t\Big|\;\oo_{\varphi}^n\geq
(1+\ga)\int_M\; \Big|\Delta_{\varphi}\pd
{\varphi}t\Big|\;\oo_{\varphi}^n,\quad t\in [\frac {\tau}2, T],$$
where $\ga>0$ is the constant obtained in Theorem
\ref{theo:prestable}. Then by the same argument in the proof of
Theorem \ref{theo:Ric1} we have the decay of the Calabi energy \beq
Ca(t)\leq e^{-\ga (t-\frac {\tau}2)}Ca(\frac {\tau}2)\leq
H^2\ee^2e^{-\ga (t-\frac {\tau}2)}, \quad \forall t\in [\frac
{\tau}2, T].\eeq Applying the parabolic Moser iteration as in the
proof of Theorem \ref{theo:Ric1}, we get the exponential decay of
the traceless Ricci curvature. The theorem is established.

\end{proof}
\subsection{The curvature tensor}
In this section, we use the exponential decay of the traceless Ricci
curvature to estimate the curvature tensor along the K\"ahler-Ricci
flow.

\begin{theo}\label{theobisec}Suppose along the K\"ahler-Ricci flow, the traceless Ricci curvature
\beq |Ric-\oo|(t)\leq H\ee e^{-\al t},\quad t\in [0, T],
\label{Ric}\eeq then the full curvature tensor is uniformly bounded
$$|Rm|(t)\leq \La_0(H\ee, \al, \La), \quad t\in [0, T],$$
where $\La$ is the bound of the curvature tensor with respect to the
background metric $\oo_g.$
\end{theo}

\begin{rem}
Phong-Song-Sturm-Weinkove proved a similar result in their paper
(cf. Lemma 6 in \cite{[PSSW]}). However, the proof here is
elementary and we don't need to use Perelman's estimates.
\end{rem}
\medskip

First we recall a useful normalization for the solution $\varphi(t)$
of the K\"ahler-Ricci flow. Observe that for any solution
$\varphi(t)$ of K\"ahler Ricci flow
$$\pd {\varphi}t=\log \frac {\oo^n_{\varphi}}{\oo_g^n}+\varphi-h_{g},$$
the function $\tilde\varphi(t)=\varphi(t)+Ce^t$ also satisfies the
above equation for any constant $C$. Since $\pd {\tilde
\varphi}{t}(0)=\pd {\varphi}{t}(0)+C,$ we have $\tilde c(0)=c(0)+C$,
where $$\td c(t)=\frac 1V\int_M\; \pd {\td \varphi}t
\oo_{\varphi}^n,\quad c(t)=\frac 1V\int_M\; \pd { \varphi}t
\oo_{\varphi}^n.$$ Thus we can normalize the solution $\varphi(t)$
such that the average of $\dot\varphi(0)$ is any given constant.\\

Recall that in Chen-Tian's paper \cite{[CT1]}, if the $K$-energy is
bounded from below, we can normalize the solution such that $c(t)$
is uniformly bounded for all $t>0$. However, in our case we have no
information about the $K$-energy. To overcome this difficulty, we
estimate $c(t)$ as follows:

\begin{lem}\label{lemchentian}(cf. \cite{[CT1]})Suppose that for $t\in [0, T]$
$$\mu_1(t)=\frac 1V\int_M\; \Big|\Na \vp\Big|^2 \oo^n_{\varphi}\leq C_4 e^{-\al t},$$
for some constant $C_4>0.$ Then we can normalize the solution
$\varphi(t)$ so that
$$c(0)=\frac 1V\int_0^{T}\;e^{-t}\int_M\;\Big|\Na \vp\Big|^2\oo_{\varphi}^n\wedge dt<C_4. $$
Then  for all time $t\in [0, T]$, we have the following estimates
$$0<c(t)<C_4 e^{-\al t},\quad \int_0^{T}\;c(t)dt<\frac {C_4}{\al}.$$
\end{lem}
\begin{proof} A simple calculation yields
$$c'(t)=c(t)-\frac 1V\int_M\;\Big|\Na \vp\Big|^2\oo_{\varphi}^n.$$
Now we normalize our initial value of $c(t)$ as \beqs c(0)=\frac
1V\int_0^{T}\;e^{-t} \int_M\;\Big|\Na \vp\Big|^2\oo_{\phi}^n\wedge
dt\leq C_4\int_0^{T}e^{-(1+\al) t} dt\leq\frac {C_4}{1+\al}\leq C_4.
\eeqs From the equation for $c(t)$, we have
$$(e^{-t}c(t))'=-\mu_1(t)e^{-t}.$$
Thus, we have \beqs0<c(t)=\frac 1V \int^{T}_t \;e^{-(\tau-t)}
\int_M\;\Big|\Na \vp\Big|^2(\tau) \oo^n_{\varphi}\wedge d\tau\leq
C_4e^{-\al t} \eeqs and
$$\int_0^{T}\;c(t)dt=C_4\int_0^{T}\; e^{-\al t} dt\leq \frac {C_4}{\al}.$$
\end{proof}
\medskip
\begin{proof}[Proof of Theorem \ref{theobisec}] We will use Theorem
\ref{theoRm} to bound the full curvature tensor. It suffices to
bound $\varphi, \pd {\varphi}t$ for time $t\in [0, T].$ Here we need
to normalize $\varphi(t)$ such that its average $c(t)$ has good
estimates. Note that the normalization in Lemma \ref{lemchentian}
depends on $T$, which is different from \cite{[CT1]}.

By the assumption (\ref{Ric}) on the Ricci curvature,  we have
$$\mu_1(t)\leq 2Ca(t)\leq 2H^2\ee^2e^{-2\al t}, \quad t\in [0, T],$$
where we used that $\la_1>1$ under the assumption of the main
theorems, or by Lemma 4.13 in \cite{[CLW]} the first eigenvalue
$\la_1(t)\geq\frac 12$ if we choose $\ee$ sufficiently small.
 By Lemma \ref{lemchentian}, we can normalize the solution
$\varphi_1(t)=\varphi(t)+G_1 e^t$, where
$$G_1=\frac 1V\int_0^{T}e^{-t}\int_M\; \Big|\Na \pd {\varphi}t\Big|^2
\oo_{\varphi}^n\wedge dt+\frac 1{V}\int_M\; h_g\; \oo_g^n,$$ such
that
$$0< c_1(t)=\frac 1V\int_M\; \pd {\varphi_1}t\;\oo_{\varphi}^n<2H^2\ee^2e^{-2\al t},
 \quad \int_0^{T}c_1(t)dt\leq \frac {H^2\ee^2}{\al}.$$
Since at the initial time $\varphi(0)=0$ and $|h_g|_{C^0}\leq
C(\sigma)H\ee$ by the normalization condition (\ref{norh}), we have
\beq |\varphi_1|(0)\leq |G_1|\leq \frac
{2H^2\ee^2}{1+2\al}+C(\sigma)H\ee.\label{sec3.4:eq1}\eeq Recall that
$$\Delta_{\varphi}\Big(\pd {\varphi_1}t-c_1(t)\Big)=n-R,$$
we have \beq \Big|\pd {\varphi_1}t-c_1(t)\Big|\leq C(\sigma)H\ee
e^{-\al t}.\label{sec3.4:eq2}\eeq Then (\ref{sec3.4:eq1}) and
(\ref{sec3.4:eq2}) imply
$$\Big|\pd {\varphi_1}t\Big|\leq\Big|\pd {\varphi_1}t-c_1(t)\Big|+|c_1(t)| \leq  2H^2\ee^2e^{-2\al t}+C(\sigma)H\ee e^{-\al t},$$
and \beqs |\varphi_1|(t)&\leq& |\varphi_1|(0)+\Big|\int_0^T\Big(\pd
{\varphi_1}t-c_1(t)\Big)dt\Big|+\int_0^T\,c_1(t)\,dt\\&\leq&2H^2\ee^2+C(\sigma)H\ee+\frac
{H^2\ee^2}{\al}+\frac {C(\sigma)}{\al}H\ee.\eeqs  Then by Theorem
\ref{theoRm} we have the curvature bound
$$|Rm|(t)\leq \La_0(H\ee, \al, \La).$$
Here we omit the Sobolev constant in the constant $\La_0$ because
the traceless Ricci curvature is very small.

\end{proof}

\section{Proof of Theorem \ref{main5}}
In this section, we follow the idea of our previous paper
\cite{[CLW]} to prove Theorem \ref{main5}. The proof needs the
technical condition that the first eigenvalue of the initial metric
is strictly greater than 1, which will be removed in Section
\ref{Sec6} by the pre-stable condition.
\medskip

\subsection{The traceless Ricci curvature}
In this subsection, we will prove that under the assumption of
Theorem \ref{main5} the traceless Ricci curvature is small after a
short time. In fact, we have the following proposition:\\

\begin{prop}\label{sec7:prop}Given a K\"ahler metric $\oo_g$ satisfying the properties
\beq \la_1(\oo_g)>1+\la,\quad |Rm|(\oo_g)\leq \La, \quad
Ca(\oo_g)<\ee, ,\label{eq:sec7:1}\eeq for some constants $\la,
\La>0$ and sufficiently small $\ee>0.$ Then for the solution
$\oo(t)$ of the K\"ahler-Ricci flow with the initial metric $\oo_g$,
there exists $T(\La)$ and $ \ee_0(\ga, \La)>0$ such that if $\ee\in
(0, \ee_0)$ then at time $T$ the metric $\oo(T)$ satisfies \beq
\la_1(T)>1+\frac {\ga}2,\quad |Rm|(T)\leq 2\La, \quad
|Ric-\oo|(T)\leq
C\ee^{\frac 14},\eeq for some constant $C(\La).$\\
\end{prop}

First, we will prove that the Sobolev constant is uniformly bounded
if the Calabi energy is small and the Ricci curvature has a lower
bound. We will follow an approach taken by \cite{[Sp]} and also
\cite{[CT2]}. Recall that C. Sprouse proved the following lemma in
\cite{[Sp]}:

\begin{lem}\label{sec7:lemSp}Let $(M, g)$ be a complete m-dimensional Riemannian manifold with $Ric\geq (m-1)k(k\leq
0).$ Then for any $D, \dd>0$ there exists $\ee= \ee(n, k, D, \dd)$
such that if \beq \sup_x \frac 1{\Vol(B(x, D))}\int_{B(x, D)}\;
((m-1)-Ric_-)_+dV<\ee,\eeq then $(M, g)$ is compact, with
$diam(M)<\pi+\dd.$ Here $Ric_-$ denotes the lowest eigenvalue of the
Ricci tensor. For any function $f$ on $M$, $f_+=\max\{f(x), 0\}.$

\end{lem}

\medskip

\begin{rem}After rescaling, the conclusion of Lemma \ref{sec7:lemSp} should be: Suppose  $Ric\geq -\La.$ For
any $D, \dd, a>0$ there exists $\ee_0=\ee_0(m, D, \dd, \La, a)>0$
such that if \beq \sup_x \frac 1{\Vol(B(x, D))}\int_{B(x, D)}\;
|Ric-a|\,dV<\ee_0,\eeq then the diameter of the metric $\oo$ is
bounded by $\sqrt{\frac {m-1}{a}}(\pi+\dd).$\\
\end{rem}

 Now we use Lemma \ref{sec7:lemSp} to give a uniform upper
bound of the Sobolev constant.

\medskip
\begin{lem}\label{sec7:lemsob}Let $(M, [\oo])$ be a polarized K\"ahler manifold and
$[\oo]$ is the canonical K\"ahler class. For any $\La>0$ there
exists $\ee_0=\ee_0(\La)$ such that if \beq Ric(\oo)\geq -\La, \quad
Ca(\oo)<\ee_0,\eeq then the Sobolev constant of the metric $\oo$ is
uniformly bounded by some constant $\si=\si(\La).$
\end{lem}

\begin{proof}Let $D, \dd$ be some fixed number, for example $D=100,
\dd=1.$ Then for any $x\in M,$ we have \beqs  \frac 1{\Vol(B(x,
D))}\int_{B(x, D)}\; |Ric-\oo|\,dV&\leq& C(n, D, \La)\Big(\int_M\;
|Ric-\oo|^2\,\oo^n\Big)^{\frac 12}\\
&\leq &C(n, D, \La)\ee_0^{\frac 12},\eeqs where we have used the
volume comparison theorem. Then by Lemma \ref{sec7:lemSp} the
diameter is bounded by a constant $C(n)$ if $\ee_0$ is small enough.
Since the volume is fixed and $Ric\geq -\La, $ by Theorem
\ref{theo:Sob} the Sobolev constant is uniformly bounded by
$\si(\La)$.\\
\end{proof}

\medskip

Now we return to the proof of Proposition \ref{sec7:prop}.

\begin{proof}[Proof of Proposition \ref{sec7:prop}] We start with the metric $\oo_g$ satisfying (\ref{eq:sec7:1}).
By the maximum principle, there exists $T(\La)>0$ such that along
the flow the curvature has the estimates \beq |Rm|(t)\leq 2\La,
\quad t\in [0, T].\label{sec7:eq:Rm}\eeq  Now we show the Calabi
energy is also small for $t\in [0, T]:$

\begin{lem}\label{sec7:lem1}Along the K\"ahler-Ricci flow, if $|Rm|(t)\leq \La, t\in [0,
T]$ then the Calabi energy satisfies
$$Ca(t)\leq e^{C(\La)t}Ca(0),\quad t\in [0, T].$$
\end{lem}
\begin{proof}By the curvature assumption and (\ref{eq:Ric1:1}) we get
\beqn \frac d{dt}Ca(t)&=&\int_M\; \Big(-2\Delta_{\varphi}\vp
\Big|\Na\bar \Na \vp\Big|^2-2\Big|\Na\Delta_{\varphi}
\vp\Big|^2+2\Big|\Delta_{\varphi} \vp\Big|^2+\Big(\Delta_{\varphi}
\vp\Big)^3\Big)\;\oo_{\varphi}^n\nonumber \\&\leq &C(\La)Ca(t).
 \eeqn
The lemma follows immediately.
\end{proof}

\medskip

 Since the initial Calabi energy is small, it follows from Lemma
 \ref{sec7:lem1} and (\ref{sec7:eq:Rm}) that
\beq Ca(t)\leq C(\La, T)\ee, \quad t\in [0, T].
\label{sec7:eq:Ca}\eeq Then by Lemma \ref{sec7:lemsob} the Sobolev
constant is uniformly bounded for time $t\in [0, T]:$  \beq
C_S(t)\leq \si(\La),\quad t\in [0, T].\eeq Combining this with the
parabolic moser iteration as in the proof of Theorem
\ref{theo:Ric1}, we have \beqn |Ric-\oo|(T)\leq C(\La, \si,
T)\Big(\int_{\frac T2}^T\; Ca(s)\;ds\Big)^{\frac 14}\leq C(\La, \si,
T)\ee^{\frac 14}. \eeqn So it remains to prove $\la_1(T)>1+\frac
{\ga}2$ provided $\ee$ is sufficiently small. We remind the readers
that it cannot be proved by Lemma \ref{lemeig} since the traceless
Ricci curvature may not be small near $t=0.$ However, we have the
lemma:

\medskip

\begin{lem}\label{sec7:lem2}Let $\la_1(t)$ be the first eigenvalue of $\oo_{\varphi}$ along the K\"ahler-Ricci
flow. Then for any constants $\La, \si>0$ there exists $\ee_0(\La,
\si)>0$ small enough  such that if for all $t\in [0, T]$ \beq
Ca(\oo_{\varphi})<\ee_0,\quad |Ric|(\oo_{\varphi})\leq \La, \quad
C_S(\oo_{\varphi})<\si,\eeq we have
$$\la_1(t)\geq \frac {\la_1(0)}{1+\la_1(0)(1-e^{-3\ee t})}e^{-3\ee t},\quad \forall t\in [0,
T],$$ where $\ee=3C(\La, n)\si\ee_0^{\frac 1n}.$
\end{lem}

\begin{proof}It follows from (\ref{eq:lemeig1}) that
 \beqn \frac {d\la_1(t)}{dt}
&=& \int_M\; (Ric(\oo_{\varphi})-\oo_{\varphi})(\Na f, \Na
f)\;\oo_{\varphi}^n-\int_M\; R|\Na f|^2
\oo_{\varphi}^n+\la_1\int_M\;R f^2\oo_{\varphi}^n .\label{eq:
eigen4}\eeqn Recall the Sobolev inequality: \beq\Big(\int_M\;
|f|^{p}\;\oo_{\varphi}^n\Big)^{\frac 2p}\leq C_S\int_M\; (|\Na
f|^2+f^2)\oo_{\varphi}^n,\label{eq:eigen2}\eeq where $p=\frac
{2n}{n-1}.$ Now we calculate : \beqs \int_M\;
(Ric(\oo_{\varphi})-\oo_{\varphi})(\Na f, \Na
f)\;\oo_{\varphi}^n&\geq & -\Big(\int_M
|Ric(\omega_\varphi)-\omega_\varphi|^n\Big)^{\frac 1n}
\cdot \Big(\int_M |\nabla f|^p\Big)^{\frac 2p}\\
&\geq & -C(\La, n) {Ca(\omega_\varphi)}^{1\over n}  \cdot
\Big(\int_M |\nabla f|^p\Big)^{\frac 2p} \eeqs
 where the
last step follows from the assumption on Ricci curvature. Now using
the Sobolev inequality (\ref{eq:eigen2}) we get: \beqs \Big(\int_M
|\nabla f|^p\Big)^{\frac 2p}& \leq & C_S \Big(\int_M |\nabla^2 f| +
\int_M |\nabla f|^2\Big)
\\ & = & C_S\Big( \int_M |f_{ij}|^2 + \int_M |f_{i\bar j}|^2 + |\nabla f|^2\Big)
\\ & = & C_S \Big( \int_M |f_{i\bar j}|^2 + \int_M  Ric(\oo_{\varphi})(\nabla f, \nabla f) +  \int_M |f_{i\bar j}|^2 +\int_M  |\nabla f|^2\Big)
\\ & = & C_S \Big( 2 \int_M |f_{i\bar j}|^2 +  \int_M  (Ric(\oo_{\varphi}) - \omega_{\varphi})(\nabla f, \nabla f)  + 2 \int_M |\nabla f|^2\Big)
\\ & \leq & 2  C_S \int_M\;\Big( (\triangle f)^2 + |\nabla f|^2\Big) + C_S C(\La, n)Ca(\omega_\varphi)^{1\over n}
\Big(\int_M\;|\Na f|^p\Big)^{\frac 2p} \eeqs
 If the  Calabi energy is sufficiently small, this implies that
\[
\Big(\int_M |\nabla f|^p\Big)^{\frac 2p} \leq  3  C_S \int_M \Big(
(\triangle f)^2 + |\nabla f|^2\Big)=3C_S(\la_1^2+\la_1).
\]
Hence, we have the estimate: \beq \int_M\;
(Ric(\oo_{\varphi})-\oo_{\varphi})(\Na f, \Na f)\;\oo_{\varphi}^n
\geq -\ee(\la_1^2+\la_1),\label{eq: eigen3}\eeq where $\ee=3C(\La,
n)\si\ee_0^{\frac 1n}.$ Similarly, we have the estimates \beq
-\int_M\; R|\Na f|^2=\int_M\; (n-R)|\Na f|^2-n\int_M\; |\Na f|^2\geq
-\ee(\la_1^2+\la_1)-n\la_1,\eeq and \beq\la_1\int_M\;
Rf^2=\la_1\int_M\; (R-n)f^2+n\la_1\geq -\ee
(\la_1+1)\la_1+n\la_1.\label{eq: eigen5}\eeq Combining (\ref{eq:
eigen3})-(\ref{eq: eigen5}) with (\ref{eq: eigen4}), we get
 \beqs \frac {d\la_1(t)}{dt}
&\geq& -3\ee (\la_1(t)+1)\la_1.\eeqs Applying the maximum principle,
we get
$$\la_1(t)\geq \frac {\la_1(0)}{1+\la_1(0)(1-e^{-3\ee t})}e^{-3\ee t}.$$

\end{proof}

\medskip

Now using Lemma \ref{sec7:lem2}, we have the estimate
$$\la_1(T)\geq \frac {\la_1(0)}{1+\la_1(0)(1-e^{-3\ee' T})}e^{-3\ee' T}, \quad \ee'=3C(\La, n)\si(C(\La, T)\ee^{\frac 14})^{\frac 1n}.$$
The first eigenvalue $\la_1(T)>1+\frac {\ga}2$ if $\ee$ is small
enough. The proposition is established.

\end{proof}
\subsection{The iteration argument  }

Thanks to Proposition \ref{sec7:prop}, we only need to prove the
convergence of K\"ahler-Ricci flow under the assumption that the
traceless Ricci curvature is sufficiently small. In other words,
Theorem \ref{main5} follows immediately from the following:

\begin{theo}\label{main2}Let $(M,  J)$ be a K\"ahler manifold with
$ c_1(M)>0.$ For any $\ga, \La>0$, there exists a small positive
constant $\ee(\ga, \La)>0$ such that for any metric $g$ in the
subspace of K\"ahler metrics \beq\{\;\oo_g\in 2\pi
c_1(M)\;|\;\la_1(\oo_g)>1+\ga,\quad |Rm|(\oo_g)\leq \La, \quad
|Ric(\oo_g)-\oo_g|\leq \ee\label{main2:eq}\},\eeq  the
K\"ahler-Ricci flow with the initial metric $\oo_g$ will converge
exponentially fast to a K\"ahler-Einstein metric.
\end{theo}

\begin{proof}
The proof consists of several parts. \medskip

 \textbf{STEP 1}. In this step we give estimates near the initial time. Consider the K\"ahler-Ricci flow
(\ref{eq:flowpotential}) with the normalization condition
(\ref{norh}) and the assumption (\ref{main2:eq}). By Lemma
\ref{lemRm}, the maximum principle implies that there exists
$T_1(\La)>0$ such that \beq |Rm|(t)\leq 2\La,\quad |Ric-\oo|(t)\leq
2\ee, \quad \forall t\in [0, T_1].\label{proof:eq3}\eeq By Lemma
\ref{lemeig}, for sufficiently small $\ee$ we have the estimate
$$\la_1(t)\geq \la_1(0)e^{-6n\ee t}\geq (1+\ga)e^{-6n\ee T_1}\geq 1+\frac {2\ga}{3},\quad \forall t\in [0, T_1]. $$
Since we have the eigenvalue estiamte, by Theorem \ref{theo:Ric1}
the traceless Ricci curvature decays exponentially \beq
|Ric-\oo|(t)\leq 2C_2(2\La, \frac 23\ga, \tau)\ee e^{-\al t}, \quad
\forall t\in [\tau , T_1]\label{proof:eq1}\eeq where $\al=\frac
{\ga}6.$ Since the traceless Ricci curvature is small for $t\in [0,
\tau]$ by (\ref{proof:eq3}), we can assume (\ref{proof:eq1}) holds
for all $t\in [0, T_1].$ Now we choose $\ee$ sufficiently small such
that
$$2C_2(2\La, \frac 23\ga, \tau)\ee<1.$$ By Theorem \ref{theobisec},
the full curvature tensor has the estimate
$$|Rm|(t)\leq \La_0(1, \frac {\ga}6),\quad t\in [0, T_1].$$
Set $\La_1:=\max\{2\La, \La_0(1, \frac {\ga}8)\}$ and
$H_1=2C_2(\La_1, \frac 23\ga, \tau)$. We choose $\ee$ small such
that $H_1\ee<\frac 12.$ By (\ref{proof:eq1}) we have
$$|Ric-\oo|(t)\leq H_1\ee e^{-\al t}, \quad
\forall t\in [0 , T_1].$$

\medskip
\medskip

\textbf{STEP 2}. Recall that in step 1 we have the estimates at time
$t=T_1,$ \beq |Rm|(T_1)\leq \La_1,\quad |Ric-\oo|(T_1)\leq H_1\ee
e^{-\al T_1},\quad \la_1(T_1)\geq 1+\frac {2\ga}3.\eeq By Lemma
\ref{lemRm}, the maximum principle implies that there exists
$T_2(\La_1)>0$ such that \beq |Rm|(t)\leq 2\La_1,\quad
|Ric-\oo|(t)\leq 2H_1\ee e^{-\al T_1},\quad t\in [T_1,
T_1+T_2].\label{sec4:eq3}\eeq By Lemma \ref{lemeig} for sufficiently
small $\ee$ the first eigenvalue has the estimate \beq \la_1(t)\geq
\la_1(T_1)e^{-3n H_1\ee T_2}\geq 1+\frac {\ga}2, \quad \forall t\in
[T_1, T_1+T_2].\eeq Then by (\ref{eq:Calabi}) the Calabi energy
decays exponentially for all time $t\in [0, T_1+T_2]:$ \beq
Ca(t)\leq e^{-\frac {\ga}2 t}Ca(0)=H_1^2\ee^2 e^{-\frac {\ga}2 t},
\quad t\in [0, T_1+T_2].\label{sec4:eq4}\eeq As in the proof of
Theorem \ref{theo:Ric1}, the parabolic Moser iteration implies the
traceless Ricci curvature \beqs |Ric-\oo|(t)&\leq& C_2(2\La_1, \frac
{\ga}2, \tau)\Big(\int_{t-\tau}^t\; \int_M\;
|Ric-\oo|^4(s)\,\oo_{\varphi}^n\wedge ds\Big)^{\frac 14}\\&\leq&
C_2(2\La_1, \frac {\ga}2, \tau)H_1\ee e^{-\al \frac {T_1}2}e^{-\frac
{\ga}8t}, \quad t\in [T_1, T_1+T_2],\eeqs where we have used the
estimates (\ref{sec4:eq3}) and (\ref{sec4:eq4}).  Set
$H_2=\max\{H_1, C_2(2\La_1, \frac {\ga}2, \tau)H_1 e^{-\al \frac
{T_1}2}\}.$ Then we get the inequality \beq |Ric-\oo|(t)\leq H_2\ee
e^{-\frac {\ga}8t}, \quad t\in [0, T_1+T_2].\eeq
 Therefore,  by Theorem \ref{theobisec} and the definition of $\La_1$, we can choose $\ee$ small
such that $H_2 \ee<1$ and \beq |Rm|(t)\leq \La_1,\quad t\in [0,
T_1+T_2].\eeq
\medskip

 \textbf{STEP 3}. Following the argument in STEP 2, we have the
 following lemma:

\begin{lem}\label{sec4:lem11}Suppose for some $k\in \NN$ the following estimates hold:
$$|Rm|(t)\leq \La_1,\quad |Ric-\oo|(t)\leq H_k \ee e^{-\frac {\ga}8 t},\quad \la_1(t)\geq 1+\frac {\ga}2, \quad t\in [0, T_1+(k-1)T_2],$$
then there exists $k_0\in \NN$  such that if $k\geq k_0$ and $\ee$
is sufficiently small (depending on $k_0$),  the above estimates
still hold for all $t\in [T_1+(k-1)T_2, T_1+kT_2].$
\end{lem}
\begin{proof}The maximum principle and the definition of $T_2$ imply that
\beq |Rm|(t)\leq 2 \La_1,\quad |Ric-\oo|(t)\leq 2 H_k \ee e^{-\frac
{\ga}8 (T_1+(k-1)T_2)}, \quad t\in [T_1+(k-1)T_2, T_1+kT_2].
\label{sec4:eq1}\eeq Then by Lemma \ref{lemeig} the first eigenvalue
\beqs \la(t)&\geq&\la(T_1+(k-1)T_2)e^{-3n H_k \ee
T_2}\\&\geq&\la(0)e^{-\frac {24n}{\ga}H_k \ee}e^{-3n H_k \ee
T_2}\\
&\geq& 1+\frac {\ga}2, \quad t\in [T_1+(k-1)T_2, T_1+kT_2],\eeqs if
we choose $\ee$ sufficiently small (depending on $H_k$). Here we
have used the inequality
$$\la(T_1+(k-1)T_2)\geq \la_1(0)e^{-\frac {24n}{\ga}H_k \ee(1-e^{-\frac {\ga}8 t})}\geq \la(0)e^{-\frac {24n}{\ga}H_k \ee}. $$
 Then by
(\ref{eq:Calabi}) the Calabi energy decays exponentially \beq
Ca(t)\leq e^{-\frac {\ga}2 t}Ca(0)=H_k^2\ee^2 e^{-\frac {\ga}2 t},
\quad t\in [0, T_1+kT_2],\label{sec4:eq2}\eeq Combining
(\ref{sec4:eq1}) with (\ref{sec4:eq2}), we can estimate the
traceless Ricci curvature by the parabolic Moser iteration as in
step 2, \beqs |Ric-\oo|(t)&\leq& C_2(2\La_1, \frac {\ga}2,
\tau)\Big(\int_{t-\tau}^t\; \int_M\;
|Ric-\oo|^4(s)\,\oo_{\varphi}^n\wedge ds\Big)^{\frac 14}\\&\leq&
H_{k+1}\ee e^{-\frac {\ga}8t}, \quad t\in  [T_1+(k-1)T_2,
T_1+kT_2],\eeqs where $ H_{k+1}$ is defined by
$$ H_{k+1}=H_k\cdot C_2(2\La_1, \frac {\ga}2, \tau)\Big(2e^{-\frac {\ga}8(T_1+(k-1)T_2)}\Big)^{\frac 12}\Big(\frac 2{\ga}e^{\frac {\ga}2\tau}\Big)^{\frac 14}.$$
Then there exists $k_0\in \NN$ such that for all $k\geq k_0$
$$C_2(2\La_1, \frac {\ga}2, \tau)\Big(2e^{-\frac {\ga}8(T_1+(k-1)T_2)}\Big)^{\frac 12}\Big(\frac 2{\ga}e^{\frac {\ga}2\tau}\Big)^{\frac 14}\leq 1.$$
Hence, we have the estimate for $k\geq k_0,$ \beq |Ric-\oo|(t)\leq
H_{k_0} \ee e^{-\frac {\ga}8 t},\quad t\in [0, T_1+kT_2].\eeq We
choose $\ee$ small such that $H_{k_0}\ee<1$. Therefore, by Theorem
\ref{theobisec} and the definition of $\La_1$ we have
$$|Rm|(t)\leq \La_1,\quad t\in [0, T_1+kT_2].$$
The lemma is proved.

\end{proof}

\medskip

\textbf{STEP 4}. By Lemma \ref{sec4:lem11}, the bisectional
curvature is uniformly bounded and the traceless Ricci curvature
decays exponentially, and therefore the $W^{1, 2}$ norm of $\pd
{\varphi}{t}-c(t)$ decays exponentially. Then following the argument
in \cite{[CT1]}, the K\"ahler-Ricci flow will converge exponentially
fast to a K\"ahler-Einstein metric. This theorem is proved.

\end{proof}

\section{Proof of Theorem \ref{main1}}
In this section, we prove Theorem \ref{main1}. The key observation
is that if a K\"ahler metric $g$ with possibly different complex
structure $J$ is sufficiently close to a K\"ahler-Einstein metric
$(g_{KE}, J_{KE})$, the conditions in Theorem \ref{main5} or Theorem
\ref{main2} are automatically satisfied, and  Theorem \ref{main1}
follows immediately.
\medskip

To verify the condition on the first eigenvalue, we need the
following result. Let $g$ be a Riemannian metric on an
$2n$-dimensional Riemannian manifold with the first eigenvalue
$\la_1(g)$. Define the space of Riemannian metrics
$$\cA_{\dd}=\{g'   \;| (1-\dd)g\leq g'\leq (1+\dd)g\}$$
for small $\dd>0.$ We have the following lemma

\medskip
\begin{lem}\label{lem:eig}For any $\ee>0$, there exists a $\dd_0(\ee)>0$ such that for any $h\in \cA_{\dd}$ with $0<\dd<\dd_0$,
 the first eigenvalue of $h$ satisfies
$$\la_1(h)>\la_1(g)(1-\ee).$$
\end{lem}
\begin{proof}
In fact, for any smooth function $f$ with $\int_M\; f dV_h=0,$ we
have \beqs \frac 1{V_{h}}\int_M\; |\Na f|_h^2dV_h&\geq& \frac
{(1-\dd)^n}{1+\dd}\frac
1{V_{h}}\int_M\; |\Na f|_g^2dV_g\\
&\geq&\la_1(g)\frac {(1-\dd)^n}{1+\dd}\frac {V_g}{V_h}\Big(\frac
1{V_{g}}\int_M\; f^2 dV_g-\Big(\frac 1{V_{g}}\int_M\; f\;
dV_g\Big)^2\Big).\eeqs Notice that \beqs \Big(\int_M\; f\;
dV_g\Big)^2&\leq &\Big(\int_M\; fdV_h\Big)^2+ \Big(\int_M\;
f\Big(\frac {dV_g}{dV_h}-1\Big)dV_h\Big)^2.\eeqs By the assumption,
$$\frac 1{(1+\dd)^n}-1\leq \frac {dV_g}{dV_h}-1\leq \frac 1{(1-\dd)^n}-1,$$
therefore, \beqs \Big(\int_M\; f\; dV_g\Big)^2&\leq
&C_{\dd}^2V_h\int_M\; f^2dV_h\eeqs where $C_{\dd}=\max\{\frac
1{(1-\dd)^n}-1, 1-\frac 1{(1+\dd)^n}\}$. Then \beqs &&\frac
1{V_{h}}\int_M\; |\Na f|_h^2dV_h\\&\geq&\la_1(g)\frac
{(1-\dd)^n}{1+\dd}\frac {V_g}{V_h}\Big(\frac 1{V_{g}}\frac
1{(1+\dd)^n}\int_M\; f^2 dV_h-\Big(\frac 1{V_{g}}\int_M\; f\;
dV_g\Big)^2\Big)\\
&\geq&\la_1(g)\frac {(1-\dd)^n}{1+\dd}\frac {V_g}{V_h}\Big(\frac
{V_h}{V_{g}}\frac
1{(1+\dd)^n}-\frac{C_{\dd}^2V_h^2}{V_g^2}\Big)\frac {1}{V_h}\int_M\;
f^2 dV_h \\
&\geq&\la_1(g)(1-\ee)\frac {1}{V_h}\int_M\; f^2 dV_h.\eeqs This
implies
$$\la_1(h)\geq \la_1(g)(1-\ee).$$

\end{proof}

\begin{proof}[Proof of Theorem \ref{main1}]
Let $(g, J)$ be a K\"ahler metric with \beq \|(g, J)-(g_{KE},
J_{KE})\|_{C^2}\leq \ee,\label{main1:eq1}\eeq for sufficiently small
$\ee>0$, we need to check that $(g, J)$ satisfies the assumption
(\ref{main2:eq}) of Theorem \ref{main2}. In fact, since $(g_{KE},
J_{KE})$ has no holomorphic vector fields, the first eigenvalue of
the Laplacian $\Delta_{KE}$ is strictly greater than 1, and by Lemma
\ref{lem:eig} $\la_1(g)\geq\la_1(g_{KE})(1-\ee)>1.$ The rest of the
assumption (\ref{main2:eq}) can be easily verified by
(\ref{main1:eq1})(cf. Lemma 2.7, 2.8 in \cite{[F]}). The theorem is
proved.

\end{proof}

\section{Proof of Theorem \ref{main3} and \ref{main6}}\label{Sec6}
In this section, we will use the pre-stable condition to drop the
assumptions that $M$ has no nonzero holomorphic vector fields, and
the dependence of the initial first eigenvalue of the Laplacian. The
idea of the proof is similar to that in section 4.\\

First, we will prove Theorem \ref{main6}. As in Proposition
\ref{sec7:prop} we prove that the traceless Ricci curvature is small
along the K\"ahler-Ricci flow after a short time.

\begin{prop}\label{sec7:prop2}Given a K\"ahler metric $\oo_g$ satisfying the properties
\beq \quad |Rm|(\oo_g)\leq \La, \quad Ca(\oo_g)<\ee, \eeq for some
constants $ \La>0$ and sufficiently small $\ee>0.$ Then for the
solution $\oo(t)$ of the K\"ahler-Ricci flow with the initial metric
$\oo_g$, there exists $T(\La)$ and $ \ee_0( \La)>0$ such that if
$\ee\in (0, \ee_0)$ then at time $T$ the metric $\oo(T)$ satisfies
\beq |Rm|(T)\leq 2\La, \quad |Ric-\oo|(T)\leq C\ee^{\frac 14},\eeq
for some constant $C(\La).$
\end{prop}

We see that this proposition follows immediately from the proof of
Proposition \ref{sec7:prop}. Hence, Theorem \ref{main6} is a direct
corollary of the following result:

\begin{theo}\label{main4}Suppose $(M, J)$ is pre-stable and the
Futaki invariant of the class $2\pi c_1(M)$ vanishes. For any
$\La>0$, there exists $\ee(\La)>0$ such that for any metric $g$ with
its K\"ahler form $\oo_g$ in the subspace of K\"ahler metrics \beq
\{\;\oo_{g}\in 2\pi c_1(M) \; |\;|Rm|(\oo_g)\leq \La, \quad
|Ric(\oo_g)-\oo_g|\leq \ee\},\label{main4:eq}\eeq the K\"ahler Ricci
flow with the initial metric $\oo_g$ will converge exponentially
fast to a K\"ahler-Einstein metric.
\end{theo}

\medskip
\begin{proof} We follow closely the proof of Theorem \ref{main2}.

\medskip

\textbf{STEP 1}. In this step we give estimates near the initial
time. Consider the K\"ahler-Ricci flow (\ref{eq:flowpotential}) with
the normalization condition (\ref{norh}) and the assumption
(\ref{main4:eq}). By Lemma \ref{lemRm}, the maximum principle
implies that there exists $T_1(\La)>0$ such that \beq |Rm|(t)\leq
2\La,\quad |Ric-\oo|(t)\leq 2\ee, \quad \forall t\in [0,
T_1].\label{proof2:eq3}\eeq By Shi's estimates(cf. \cite{[Shi]}),
all the derivatives of the curvature tensor are uniformly bounded
for time $t\in [\frac {\tau}2, T_1].$ Then applying  Theorem
\ref{theo:prestable}, if $\ee$ is sufficiently small there exists
$\ga_0:=\ga(2\ee, 2\La)>0$ for time $t\in [\frac {\tau}2, T_1]$ such
that the inequalities (\ref{theo:pre eq1}) and (\ref{theo:pre eq2})
hold. Since we have the eigenvalue estimate, by Theorem
\ref{theo:Ric2} the traceless Ricci curvature decays exponentially
\beq |Ric-\oo|(t)\leq 2C_3(2\La, \frac {\ga_0}4, \tau)\ee e^{-\frac
{\ga_0}4t}, \quad \forall t\in [\tau, T_1].\label{proof2:eq1}\eeq
Since the traceless Ricci curvature is small for $t\in [0, \tau]$ by
(\ref{proof2:eq3}), we can assume (\ref{proof2:eq1}) holds and also
$\ga\geq \ga_0$ for all $t\in [0, T_1].$ Now we choose $\ee$
sufficiently small such that
$$2C_3(2\La, \frac {\ga_0}4, \tau)\ee<1.$$ By Theorem \ref{theobisec},
the full curvature tensor has the estimate
$$|Rm|(t)\leq \La_0(1, \frac {\ga_0}4),\quad t\in [0, T_1].$$
Set $\La_1:=\max\{2\La, \La_0(1, \frac {\ga_0)}8)\}$ and
$H_1=2C_3(2\La, \frac {\ga_0}4, \tau)$. We choose $\ee$ small such
that $H_1\ee<\frac 12.$ By (\ref{proof2:eq1}) we have
$$|Ric-\oo|(t)\leq H_1\ee e^{-\frac {\ga_0}4 t}, \quad
\forall t\in [0 , T_1].$$

\medskip
\medskip

\textbf{STEP 2}. Recall that in step 1 we have the estimates at time
$t=T_1,$ \beq |Rm|(T_1)\leq \La_1,\quad |Ric-\oo|(T_1)\leq H_1\ee
e^{-\frac {\ga_0}4 T_1}.\eeq By Lemma \ref{lemRm}, the maximum
principle implies that there exists $T_2(\La_1)>0$ such that
$$|Rm|(t)\leq 2\La_1,\quad |Ric-\oo|(t)\leq 2H_1\ee e^{-\frac {\ga_0}4 T_1},\quad t\in [T_1, T_1+T_2].$$
By Theorem \ref{theo:prestable}, we can choose $\ee$ sufficiently
small such that $\ga(2H_1\ee, 2\La_1)\geq \ga_1$ for some constant
$\ga_1\in (0, \ga_0).$ Then by (\ref{eq:Calabi}) the Calabi energy
decays exponentially for all time $t\in [0, T_1+T_2]:$ \beq
Ca(t)\leq e^{-\ga_1 t}Ca(0)=H_1^2\ee^2 e^{-\ga_1 t}, \quad t\in [0,
T_1+T_2].\eeq As in the proof of Theorem \ref{theo:Ric1}, the
parabolic Moser iteration implies the traceless Ricci curvature
\beqs |Ric-\oo|(t)&\leq& C_3(2\La_1,\ga_1,
\tau)\Big(\int_{t-\tau}^t\; \int_M\;
|Ric-\oo|^4(s)\,\oo_{\varphi}^n\wedge ds\Big)^{\frac 14}\\&\leq&
C_3(2\La_1, \ga_1, \tau)H_1\ee e^{- \frac {\ga_0T_1}8}e^{-\frac
{\ga_1}4t}, \quad t\in [T_1, T_1+T_2].\eeqs Set $H_2=\max\{H_1,
C_3(2\La_1, \ga_1, \tau)H_1e^{- \frac {\ga_0T_1}8}\},$ we get the
inequality \beq |Ric-\oo|(t)\leq H_2\ee e^{-\frac {\ga_1}4t}, \quad
t\in [0, T_1+T_2].\eeq
 Therefore,  by Theorem \ref{theobisec} and the definition of $\La_1$, we can choose $\ee$ small
such that $H_2 \ee<1$ and \beq |Rm|(t)\leq \La_1,\quad t\in [0,
T_1+T_2].\eeq
\medskip

 \textbf{STEP 3}. Following the argument in STEP 2, we have the
 following lemma:

\begin{lem}\label{sec4:lem1}Suppose for some $k\in \NN$ the following estimates hold:
$$|Rm|(t)\leq \La_1,\quad |Ric-\oo|(t)\leq H_k \ee e^{-\frac {\ga_1}4 t},\quad \forall \;t\in [0, T_1+(k-1)T_2],$$
then there exists $k_0\in \NN$  such that if $k\geq k_0$ and $\ee$
is sufficiently small (depending on $k_0$),  the above estimates
still hold for all $t\in [T_1+(k-1)T_2, T_1+kT_2].$
\end{lem}
\begin{proof}The maximum principle and the definition of $T_2$ imply that
\beq |Rm|(t)\leq 2 \La_1,\quad |Ric-\oo|(t)\leq 2 H_k \ee e^{-\frac
{\ga_1}4 (T_1+(k-1)T_2)}, \quad t\in [T_1+(k-1)T_2, T_1+kT_2].
\label{sec6:eq1}\eeq Then   if $\ee$ is sufficiently small
(depending on $H_k$), the constant $\ga$ in Theorem
\ref{theo:prestable} satisfies $$\ga(2\La_1, 2H_k)\geq \ga_1>0,
\quad t\in [T_1+(k-1)T_2, T_1+kT_2].$$
 Then by
(\ref{eq:Calabi}) the Calabi energy decays exponentially \beq
Ca(t)\leq e^{-\ga_1 t}Ca(0)=H_k^2\ee^2 e^{-\ga_1\; t}, \quad t\in
[0, T_1+kT_2].\label{sec6:eq2}\eeq Combining (\ref{sec6:eq1}) with
(\ref{sec6:eq2}), we can estimate the traceless Ricci curvature by
the parabolic Moser iteration as in step 2, \beqs |Ric-\oo|(t)&\leq&
C_3(2\La_1, \ga_1, \tau)\Big(\int_{t-\tau}^t\; \int_M\;
|Ric-\oo|^4(s)\,\oo_{\varphi}^n\wedge ds\Big)^{\frac 14}\\&\leq&
H_{k+1}\ee e^{-\frac {\ga_1}4t}, \quad t\in  [T_1+(k-1)T_2,
T_1+kT_2],\eeqs where $ H_{k+1}$ is defined by
$$ H_{k+1}=H_k\cdot C_3(2\La_1, \ga_1, \tau)\Big(2e^{-\frac {\ga_1}4(T_1+(k-1)T_2)}\Big)^{\frac 12}\Big(\frac 1{\ga_1}e^{\ga_1\tau}\Big)^{\frac 14}.$$
Then there exists $k_0\in \NN$ such that for all $k\geq k_0$
$$C_3(2\La_1, \ga_1, \tau)\Big(2e^{-\frac {\ga_1}8(T_1+(k-1)T_2)}\Big)^{\frac 12}\Big(\frac 1{\ga_1}e^{\ga_1\tau}\Big)^{\frac 14}\leq 1.$$
Hence, we have the estimate for $k\geq k_0,$ \beq |Ric-\oo|(t)\leq
H_{k_0} \ee e^{-\frac {\ga_1}4 t},\quad t\in [0, T_1+kT_2].\eeq We
choose $\ee$ small such that $H_{k_0}\ee<1$. Therefore, by Theorem
\ref{theobisec} and the definition of $\La_1$ we have
$$|Rm|(t)\leq \La_1,\quad t\in [0, T_1+kT_2].$$
The lemma is proved.

\end{proof}

\medskip

\textbf{STEP 4}. By Lemma \ref{sec4:lem1}, the bisectional curvature
is uniformly bounded and the traceless Ricci curvature decays
exponentially, and therefore the $W^{1, 2}$ norm of $\pd
{\varphi}{t}-c(t)$ decays exponentially. Then following the argument
in \cite{[CT1]}, the K\"ahler-Ricci flow will converge exponentially
fast to a K\"ahler-Einstein metric. This theorem is proved.

\end{proof}

\begin{proof}[Proof of Theorem \ref{main3}] Theorem \ref{main3} is a
direct corollary of Theorem \ref{main4}.
\end{proof}

\end{document}